\documentstyle[12pt,amscd]{amsart}
\newtheorem{pr}{Proposition}
\newtheorem{lm}{Lemma}

\newcommand{\proj}{\Bbb {P}}

\newcommand{\barr}{\overline}

\newcommand{\rarr}{\rightarrow}
\newcommand{\oh}{{\cal{O}}}
\newcommand{\M}{\barr{M}_{0,n}}
\newcommand{\MM}{\barr{M}}
\newcommand{\eqq}{\stackrel{\sim}{=}}

\newcommand{\QQ}{{\Bbb Q}}
\newcommand{\LL}{{\cal{L}}}
\newcommand{\NN}{{\cal{N}}}
\newcommand{\g}{\gamma}
\newcommand{\s}{\sigma}
\begin{document}
\title{The Symmetric function $h^0(\M, \LL_1^{x_1}\otimes
\LL_2^{x_2} \otimes \cdots \otimes \LL_n^{x_n})$  }
\author{Rahul Pandharipande$^1$}
\date{15 November 1995}
\maketitle
\pagestyle{plain}
\footnotetext[1]{Research partially 
supported by an NSF Post-Doctoral Fellowship.}
\setcounter{section}{-1}
\section{Introduction}
\subsection{Summary of Results}
Let $n\geq 3$ be an integer. Let $(C,\  p_1,\ldots, p_n)$ be a
reduced, connected, (at worst) nodal
curve of arithmetic genus $0$ with $n$ nonsingular
marked points. $(C, \ p_1, \ldots, p_n)$ is {\em Deligne-Mumford
stable} if $w_C(p_1+ \ldots +p_n)$ is ample (where $w_C$ is
the dualizing sheaf).  
Let $\M$ denote the fine moduli space of Deligne-Mumford stable,
$n$-pointed, genus $0$ curves. 
A foundational treatment of $\M$ can be found in [Kn].
Let $S_n=\{1,2, \ldots, n\}$ be 
the marking set. 
For each $i\in S_n$, a line bundle $\LL_i$ 
on $\M$ is obtained by the following prescription: the
fiber of $\LL_i$ at the moduli point 
$[C, p_1, \ldots, p_n]\in \M$ is
the cotangent space $T^*_{C,p_i}$.
\begin{pr}
\label{vanish}
Let $x_1, x_2, \ldots, x_n$ be non-negative integers.
$$\forall k>0, \ \ H^k(\M,  \LL_1^{x_1}\otimes
\LL_2^{x_2} \otimes \cdots \otimes \LL_n^{x_n})=0.$$
\end{pr}
For non-negative $x_i$, 
let $\g_n(x_1,x_2, \ldots, x_n)=h^0(\M,  \LL_1^{x_1}\otimes
\LL_2^{x_2} \otimes \cdots \otimes \LL_n^{x_n})$.
By Proposition (\ref{vanish}) and the Hirzebruch-Riemann-Roch
Theorem, $\g_n$ is a polynomial function: $\g_n \in \QQ[x_1,x_2,
\ldots, x_n]$. The symmetric group $\Sigma_n$ acts 
naturally on $\M$
by permuting the markings. The action of $\Sigma_n$
permutes the isomorphism classes of the line bundles $\LL_i$
in the obvious manner. Therefore, $\g_n$ is a {\em symmetric}
polynomial in $x_1, x_2, \ldots, x_n$. Since $\M$ is
$n-3$ dimensional, $\g_n$ is of degree (at most) $n-3$.

Let $\sigma_1, \sigma_2, \ldots, \sigma_n$ denote the
elementary symmetric functions in $x_1, x_2, \ldots, x_n$.
\begin{eqnarray*}
\sigma_1 &=& x_1+ \ldots + x_n \\
\sigma_2 &=& x_1x_2 + x_1x_3 + \ldots +x_{n-1} x_n \\
& \vdots & \\
\sigma_n & =&  x_1x_2 \cdots x_n. \\
\end{eqnarray*}
The following equivalence is well-known:
$$ \QQ[\sigma_1, \ldots, \sigma_n] = \Big( \QQ[x_1, \ldots, x_n] 
\Big)^{\Sigma_n}.$$ Hence, $\g_n\in \QQ[\sigma_1, \ldots, \sigma_n]$.

Let $R= \QQ[\sigma_1, \sigma_2, \ldots]$ be the
infinite polynomial ring in the variables $\{\sigma_d\}|_{d\in {\Bbb N}^+}$.
Define a grading on $R$ by assigning to the 
variable $\sigma_d$ the weight $d$.
For $f\in R$, the {\em degree} of $f$ is the highest weight of
the monomials that appear in $f$.

A $\QQ$-linear transformation $T: R\rarr R$ is defined
as follows. Let
$f\in R$. Let $e=degree(f)$.
Therefore, $f\in \QQ[\sigma_1, \sigma_2, \ldots, \sigma_e]$.
The element $f$ corresponds to a symmetric function in 
the variables $x_1, \ldots, x_m$ provided $m\geq e$.
Let $g(x_1, \ldots, x_m)$ be defined
for non-negative integers $x_i$ by
\begin{equation}
\label{Tdef}
g(x_1, \ldots, x_m)= f(x_1, \ldots, x_m)+
\sum _{i=1}^{m} \sum_{j=0}^{x_i-1} 
f(x_1, \ldots,x_{i-1}, j, x_{i+1}, 
\ldots, x_{m}).
\end{equation}
The function $g$ is manifestly a symmetric polynomial
in $x_1, \ldots, x_m$ of degree (at most) $e+1$. Therefore
$g\in \QQ[\sigma_1, \sigma_2, \ldots, \sigma_e, \sigma_{e+1}]$.
The element $g$ is independent
of $m$ provided $m\geq e+1$. Define $T(f)=g \in R$ for $m$
in the stable range $m\geq e+1$.

\begin{pr}
\label{solution}
The
function
$\g_n \in \QQ[\sigma_1, \ldots, \sigma_n]$ is determined by:
$$\g_n= T^{n-3}(1).$$
\end{pr}
\noindent
The author has benefited from conversations with W. Fulton,
M. Kapranov, and M. Thaddeus.

\subsection{Computing $T$}
Computing $T$ is tedious but straightforward.
If $f=1$, then $degree(f)=0$. $T(1)$ can be computed
in the variable $x_1$. By equation (\ref{Tdef}),
$$g(x_1)= f(x_1)+ \sum_{j=0}^{x_1-1} f(j) = 1+x_1 = 1+\sigma_1.$$
Hence $T(1)=1+\sigma_1 $. 

Let $n=4$. $\barr{M}_{0,4} \eqq \proj^1$. Via this
isomorphism,
$\LL_i \eqq \oh_{\proj^1}(1)$ for $1\leq i \leq 4$. The
function $\g_4$ is easily evaluated:
\begin{eqnarray*}
\g_4(x_1,x_2,x_3,x_4) & = & h^0(\proj^1, \oh_{\proj^1}(x_1+x_2+x_3+x_4)) \\
 & = &
1+x_1+x_2+x_3+x_4 \\ & = & 1+\sigma_1  \\ & = & T(1).
\end{eqnarray*}

Next, let $f=\sigma_1$.
$T(\sigma_1)$ can be computed in the variables $x_1, x_2$.
\begin{eqnarray*}
g(x_1,x_2)& =&  f(x_1,x_2)+ \sum_{j=0}^{x_1-1} f(j,x_2) 
                        + \sum_{j=0}^{x_2-1} f(x_1,j) \\
& = & x_1+x_2 +{1\over 2} (x_1^2-x_1)  + x_1x_2 + 
               {1\over 2} (x_2^2-x_2)  + x_1x_2  \\
& = & {1\over 2}\sigma_1 +{1\over 2} \sigma_1^2 +\sigma_2  \\
\end{eqnarray*}
Therefore, $T(\s_1)={1\over 2}\sigma_1 +{1\over 2} \sigma_1^2 +\sigma_2$.
By Proposition (2),
$$\g_5=T^2(1)=T(\sigma_1+1)=T(\sigma_1)+T(1)=
1+{3\over 2}\sigma_1+ {1\over 2}\sigma_1^2+\s_2.$$

Proposition (2) leads to an easy calculation of $h^0(\overline{M}_{0,n}
, \LL_1^{x_1})$. By relation (\ref{Tdef}), for $n\geq 4$,
$$h^0(\overline{M}_{0,n}
, \LL_1^{x_1}) = \sum_{j=0}^{x_1} h^0(\overline{M}_{0,{n-1}},
\LL_1^{x_1}).$$
It is easily seen $h^0(\overline{M}_{0,n}, \LL_1^{x_1})={n-3+x_1\choose
x_1}$ uniquely satisfies the above recursion and the boundary
conditions at $n=3$. The global sections of $\LL_1^{x_1}$ 
can also be computed by examining the 
morphism $\overline{M}_{0,n}\rarr \proj^{n-3}$
determined by $H^0(\overline{M}_{0,n}, \LL_1)$ (see
section (\ref{express})).

The function $\g_n$ is tabulated below for small values of $n$:

\begin{eqnarray*}
\g_3 &= &  1\\
\g_4 &=&  1+\sigma_1\\
& & \\
\g_5 &=&  1+{3\over 2}\sigma_1+ {1\over 2}\sigma_1^2+\s_2\\
& & \\
\g_6 &=& 1+{11\over 6}\s_1+\s_1^2+\s_2+{1\over 6}\s_1^3+\s_1\s_2+2\s_3\\
& & \\
& & \\
\g_7 &=& 1+{25\over 12}\s_1+{35\over 24}\s_1^2 + {5\over 4} \s_2
          +{5\over 12}\s_1^3+ {5\over 4} \s_1\s_2 -{1\over 4}\s_3 \\
     & &\ \ \     +{1\over 24}\s_1^4+{1\over 2}\s_1^2\s_2 +{5\over 2}\s_1\s_3
          +{1\over 4}\s_2^2 +{11\over 2}\s_4 \\
& & \\
& & \\
\g_8 &=& 1+{137\over 60}\s_1+ {15\over 8}\s_1^2+ {5\over 4}\s_2
          +{17\over 24}\s_1^3+{11\over 6}\s_1\s_2 +{11\over 4}\s_3 \\
& & \ \ \ +{1\over 8}\s_1^4+{3\over 4}\s_1^2\s_2 -{3\over 4}\s_1\s_3
          +{1\over 4}\s_2^2 -{21\over 2}\s_4 \\
 & & \ \ \ +{1\over 120}\s_1^5 +{1\over 6}\s_1^3\s_2 +{3\over2}\s_1^2\s_3
     +{1\over 4}\s_1\s_2^2 +{17\over 2}\s_1\s_4+\s_2\s_3 +19 \s_5 \\
\end{eqnarray*}

\section{The Line bundles $\LL_i$}
\subsection{Expressions for $\LL_i$}
\label{express}
Let $\rho:U_n\rarr \M$ be the universal family of
pointed, stable curves. The are $n$ sections
$s_1, \ldots, s_n$ of $\rho$ corresponding to the markings.
Consider the marking $i$. The projection $\rho$
takes $s_i$ isomorphically to $\M$. 
Let $$x=[C, p_1, \ldots, p_n]\in \M.$$
Let $y\in s_i$ satisfy $\rho(y)=x$.
The normal bundle to $s_i$ in $U_n$ at $y$ is canonically the
the tangent space $T_{C,p_i}$. Therefore
\begin{equation}
\label{bay}
c_1(\LL_i)\eqq \rho_*(-s_i^2).
\end{equation}

Note $U_n$ is canonically isomorphic to $\barr{M}_{0,n+1}$.
$\LL_{n+1}$ on $U_n$ can be expressed as follows:
\begin{equation}
\label{om}
\LL_{n+1}\eqq \omega_{\rho}(s_1+\ldots + s_n)
\end{equation}
where $\omega_{\rho}$ is the relative dualizing
sheaf of the morphism $\rho: U_n\eqq \barr{M}_{0,n+1}
 \rarr \M$.
The proof of (\ref{om}) requires a square diagram:
\begin{equation*}
\begin{CD}
\barr{M}_{0,n+1\#} @>{\pi_{n+1}}>> \barr{M}_{0,n+1}\\
@V{\nu}VV @V{\rho}VV \\
\barr{M}_{0,n\#} @>{\pi_{n}}>> \M 
\end{CD}
\end{equation*}
$\MM_{0,n+1\#}$ is the moduli
space of $n+2$-pointed genus 0 stable curves with marking
set $\{1, \ldots, n+1, \#\}$. 
The morphisms $\nu$, $\rho$, $\pi_n$, and $\pi_{n+1}$ are 
all contraction
maps.
Let $D_{n+1, \#}$
be the 
boundary divisor of $\MM_{0,n+1\#}$
corresponding to the partition $$\{n+1,\#\} \cup \{1, \ldots, n\}.$$
$D_{n+1, \#}$ is the section of ${\pi_{n+1}}$ corresponding to
the marking $n+1$.
A simple examination of the blow-ups involved in the
above square yields:
$$\pi_{n+1}^*(\omega_{\rho}(s_1+\ldots + s_n))|_{D_{n+1,\#}} \eqq 
\omega_{\nu}|_{D_{n+1, \#}}.$$
Since $\omega_{\nu}|_{D_{n+1, \#}} \eqq -D_{n+1, \#}|_{D_{n+1, \#}}$,
(\ref{om}) is established by (\ref{bay}).

Another method of viewing $\LL$ is as follows.
Fix $n$ marked points in general linear position in $\proj^{n-2}$.
M. Kapranov has shown the closure in the Hilbert scheme
of the locus of rational normal curves passing through
the $n$ marked points is canonically $\M$ ([K]).
The universal curve over the Hilbert scheme restricts
to $U_n$ over $\M$. Since the universal
curve over the Hilbert scheme naturally maps
to $\proj^{n-2}$, a morphism $$\mu:U_n \rarr \proj^{n-2}$$ is
obtained.
$U_n$ is canonically isomorphic to $\barr{M}_{0,n+1}$.
There is an isomorphism:
\begin{equation}
\label{pullb}
\LL_{n+1} \eqq \mu^*(\oh_{\proj^{n-2}}(1)).
\end{equation}
On each fiber of $\rho: U_n \rarr \M$,
$\mu^*(\oh_{\proj^{n-2}}(1))$ is isomorphic to
$\omega_{\rho}(s_1+\ldots + s_n)$ (see [K]). Since 
both line bundles
$\mu^*(\oh_{\proj^{n-2}}(1))$ and $\omega_{\rho}(s_1+\ldots + s_n)$
are trivial on the sections $s_i$ of $\rho$,
$$\mu^*(\oh_{\proj^{n-2}}(1)) \eqq \omega_{\rho}(s_1+\ldots + s_n)$$
on $U_n$. Now (\ref{pullb}) follows from (\ref{om}).

Finally, it is useful to express $\LL_{n+1}$ on $\MM_{0,n+1}$
as a linear combination of boundary divisors.
Let $j,k\in S_{n}$ be distinct markings.
\begin{equation}
\label{lincom}
\LL_{n+1} \eqq \sum_{n+1 \in A\subset S_{n+1},\  j,k\notin A} D_A
\end{equation}
where the set $A$ is a subset of the marking set $S_{n+1}$
satisfying $|A|\geq 2$.
As above, $D_A$ is the boundary divisor
corresponding to the partition $S_{n+1}=A\cup A^c$.
Recall the map $$\mu:\MM_{0,n+1} \eqq U_n \rarr \proj^{n-2}.$$
Consider the unique hyperplane  $H\subset \proj^{n-2}$ 
passing through
the $n-2$ marked points $S_n \setminus \{j,k\}$. 
The divisor $\mu^*(H)$ is easily seen to be the
right side of (\ref{lincom}). By (\ref{pullb}),
the isomorphism of (\ref{lincom}) is established.
See [W] for another proof of (\ref{lincom}).
\label{lineb}

\subsection{Contraction}
It will be helpful in the sequel to denote
$\LL_i$ on $\M$ by $\LL_{i,n}$.
Let $n\geq 3$.
Consider the contraction morphism
$$\rho: \MM_{0, n+1} \rarr \M$$
obtained by omitting the marking $n+1$.
\begin{lm}
\label{comp}
$\rho^*(\LL_{i,n}) + D_{i,n+1} \eqq \LL_{i, n+1}. $
\end{lm}
\begin{pf}
Let $j,k\in S_n$ satisfy $j,k\neq i$.
By the results of section (\ref{lineb})
$$\rho^*(\LL_{i,n}) \eqq \sum_{i \in A\subset S_n,\  j,k\notin A} 
\rho^*(D_A).$$
A comparison with
$$\LL_{i,n+1} \eqq \sum_{i \in A\subset S_{n+1},\  j,k\notin A} D_A$$
yields Lemma (\ref{comp}).
\end{pf}

\section{The Proof of Proposition (\ref{vanish})}
\subsection{The Induction Ladder}
Consider the following (infinite) commutative diagram:
\begin{equation}
\label{ladder}
\begin{CD}
\MM_{0,n+1\#} @>{\nu_n}>>
\barr{M}_{0,n\#} @>{\nu_{n-1}}>> \barr{M}_{0,n-1 \#}\\
@V{\pi_{n+1}}VV  @V{\pi_{n}}VV @V{\pi_{n-1}}VV \\
\MM_{0,n+1} @>{\rho_{n}}>>
\barr{M}_{0,n} @>{\rho_{n-1}}>> \MM_{0,n-1} 
\end{CD}
\end{equation}
The maps are all contraction morphisms:
$\pi_n$ contracts $\#$, $\nu_{n-1}$ contracts $n$,
$\rho_{n-1}$ contracts $n$.
The diagram starts with $\MM_{0,3}$ on the
lower right corner and extends left.
Denote the composition $\rho_n \circ \cdots \circ \rho_{m-1}$
by $\rho_{m,n}$ (similarly for $\nu$).    

Recall $\LL_{i,n}$ is the line bundle
corresponding to the $i^{th}$ marking on $\M$.
For $m>n$, the pull-back line bundle on $\MM_{0,m}$,
$\rho_{m,n}^*(\LL_{i,n})$,
is also denoted by the same symbol $\LL_{i,n}$.

Let $\LL_{\#,n\#}$ be the line bundle corresponding
to the marking $\#$ on $\MM_{0,n\#}$. Again
$\LL_{\#,n\#}$ will also denote the pull-back line
bundle, $\nu_{m,n}^*(\LL_{\#,n\#})$ on $\MM_{0,m\#}$ for $m>n$.

Proposition (1) is established by an induction on the
ladder (\ref{ladder}). The following Lemmas are
needed in the induction.
\begin{lm}
\label{d} Let $n\geq 3$. The line bundle $\LL_{\#,n\#}$ behaves well
under ladder base changes:
\begin{enumerate}
\item[(i.)]
For all $k\geq 0$ and  $m\geq n$, $\ \pi_{m*}(\LL^k_{\#,n\#})\eqq 
\rho_{m,n}^*\pi_{n*}(\LL^k_{\#,n\#}).$
\item[(ii.)]
For all $k\geq 0$ and $m\geq n$,
$\ R^1\pi_{m*} (\LL^k_{\#,n\#})=0.$
\end{enumerate}
\end{lm}
\begin{pf}
By (\ref{pullb}), $\LL_{\#,n\#}$ is generated by global
sections on $\MM_{0,n\#}$. Therefore $\LL^k_{\#,n\#}$ is
generated globally on
$\MM_{0,m\#}$. Any line bundle generated by global
sections on a genus 0, stable, pointed curve has
no higher cohomology. Hence, $\LL^k_{\#,n\#}$ has no
higher cohomology on the fibers of $\pi_m$.
By the Cohomology and Base Change Theorems,
$$R^1\pi_{m*} (\LL^k_{\#,n\#})=0.$$
Consider the fiber product:
$$\tau_m: \MM_{0,m}\times_ {\MM_{0,n}} \MM_{0,n\#} \rarr \MM_{0,m}.$$
Again by Base Change, $\tau_{m*}(\LL^k_{\#,n\#}) \eqq \rho_{m,n}^*
\pi_{n*}(\LL^k_{\#,n\#})$. There is 
natural map
\begin{equation}
\label{nat}
\MM_{0,m\#} \rarr \MM_{0,m}\times_ {\MM_{0,n}} \MM_{0,n\#}
\end{equation}
which commutes with $\pi_m$ and $\tau_m$.
It is easy to check that the natural map of
vector bundles $\tau_{m*}(\LL^k_{\#,n\#}) 
\rarr \pi_{m*}(\LL^k_{\#,n\#})$
induced by (\ref{nat}) is an isomorphism on fibers.
\end{pf} 

\begin{lm}
\label{dd}
Let $n\geq 4$.
The following exact sequences of
vector bundles exist on $\MM_{0,n}$ for
each pair of integers $a,b \geq 0$:
\begin{equation}
\label{w}
0 \rarr \pi_{n*}(\LL^{a+1}_{\#,n-1\#}\otimes \LL^b_{\#,n\#})
\rarr \pi_{n*}(\LL^{a}_{\#,n-1\#}\otimes \LL^{b+1}_{\#,n\#})
\rarr \LL^a_{n,n} \rarr 0.
\end{equation}
\end{lm}
\begin{pf}
By Lemma (\ref{comp}), 
there is a linear equivalence $\LL_{\#,n-1\#} +
D_{\#,n} \eqq \LL_{\#, n\#}$. 
Note: $$\LL_{\#, n\#}|_{D_{\#,n}}=0,$$
\begin{equation}
\label{jonny}
\LL_{\#,n-1\#}|_{D_{\#,n}}  \eqq -D_{\#,n}|_{D_{\#,n}} .
\end{equation}
Tensoring the 
sequence (\ref{ww}) with $\LL^a_{\#, n-1\#}\otimes 
\LL^{b+1}_{\#,n\#}$
yields (\ref{www}).
\begin{equation}
\label{ww}
0 \rarr \oh(-D_{\#,n}) \rarr \oh \rarr \oh_{D_{\#,n}} \rarr 0
\end{equation}
\begin{equation}
\label{www}
0 \rarr \LL^{a+1}_{\#,n-1\#}\otimes \LL^b_{\#,n\#}
\rarr \LL^{a}_{\#,n-1\#}\otimes \LL^{b+1}_{\#,n\#}
\rarr  \LL^a_{\#,n-1\#}|_{D_{\#,n}} \rarr 0.
\end{equation}
Equivalence (\ref{jonny}) implies
$\pi_{n*}(\LL^a_{\#,n-1\#}|_{D_{\#,n}}) \eqq \LL^a_{n,n}$.
By a global sections argument as in Lemma (\ref{d}),
the terms of sequence (\ref{www}) have vanishing
higher direct images under $\pi_n$.
Sequence
(\ref{w}) is obtained by pushing-forward sequence (\ref{www}).
\end{pf}

\begin{lm}
\label{ddd} Let $n\geq 3$.
Let $\NN$ be a line bundle on $\M$ satisfying the
following condition 
for all non-negative integers $z_3,\ldots, z_n$
and all $k>0$ :
\begin{equation}
\label{f}
H^k(\M, \NN\otimes \bigotimes_{i=3}^n \LL^{z_i}_{i,i})=0.
\end{equation}
Then for all  $b\geq 0$ and $k>0$,
\begin{equation}
\label{ff}
H^k(\M, \pi_{n*}(\LL^b_{\#,n\#})\otimes \NN)=0.
\end{equation}
\end{lm} 
\begin{pf}
If $b=0$, the vanishing (\ref{ff}) is a consequence
of (\ref{f}) since $\pi_{n*}(\oh)\eqq \oh$. Assume $b>0$.
The proof is a simple consequence of Lemmas (\ref{d}) and (\ref{dd}).
Consider
the sequences for $0\leq j \leq b-1$ obtained from (\ref{w}) by tensoring
with $\NN$:
$$0 \rarr \pi_{n*}(\LL^{j+1}_{\#,n-1\#}\otimes 
\LL^{b-j-1}_{\#,n\#})\otimes \NN
\rarr \pi_{n*}(\LL^{j}_{\#,n-1\#}\otimes 
\LL^{b-j}_{\#,n\#})\otimes \NN
\rarr \LL^j_{n,n} \otimes \NN \rarr 0.$$
The vanishing (\ref{ff}) is reduced by these sequences and
repeated application of (\ref{f}) to:
  $$\forall b\geq 0, \ \forall k>0, \ \ \
 H^k(\M, \pi_{n*}(\LL^b_{\#,n-1\#})\otimes \NN)=0.$$
But since, $\pi_{n*}(\LL^b_{\#,n-1\#})\eqq \rho_{n-1}^*
\pi_{n-1*}(\LL^b_{\#,n-1\#})$ by Lemma (\ref{d}),
it suffices to show:
$$\forall b\geq 0, \ \forall k>0, \ \ \
  H^k(\M, \rho^*_{n-1}\pi_{n-1*}(\LL^b_{\#,n-1\#})\otimes \NN)=0.$$
Pulling-back the sequences (\ref{w}) for $n-1$ to $\M$ via $\rho_{n-1}$,
tensoring with $\NN$, and repeatedly applying (\ref{f}),
the vanishing (\ref{ff}) is further reduced to:
$$\forall b\geq 0, \ \forall k>0, \ \ \
  H^k(\M, \rho^*_{n-1}\pi_{n-1*}(\LL^b_{\#,n-2\#})\otimes \NN)=0.$$
Applying Lemma (\ref{d}), it suffices to show
$$\forall b\geq 0, \ \forall k>0, \ \ \
  H^k(\M, \rho^*_{n,n-2}\pi_{n-2*}(\LL^b_{\#,n-2\#})\otimes \NN)=0.$$
This process finally reduces the original claim (\ref{ff})
to:
$$\forall b\geq 0, \ \forall k>0, \ \ \
  H^k(\M, \rho^*_{n,3}\pi_{3*}(\LL^b_{\#,3\#})\otimes \NN)=0.$$
Since $\pi_{3*}(\LL^b_{\#,3\#})$ is a trivial bundle over
the point $\MM_{0,3}$, the
Lemma is proven.
\end{pf}

\subsection{The Induction}
A slightly stronger version of Proposition (\ref{vanish}) is
needed for the induction:

\noindent {$\bold {Proposition\ 1'}.$}
Let $n\geq 3$.
Let $x_{i,j}$ (for $1\leq i \leq j$, $3\leq j \leq n$ ) 
be non-negative integers.
$$\forall k>0, \ \ H^k(\M,\ \bigotimes_{i,j} \LL_{i,j}^{x_{i,j}})=0.$$

\begin{pf}
The Proposition certainly holds for $n=3$ since $dim(\MM_{0,3})=0$.
Assume the Proposition holds for $n$.
Let $x_{i,j}$ (for $1\leq i \leq j$, $3\leq j \leq n$) be non-negative
integers. Let $y_1, \ldots, y_n, y_{\#}$ also be non-negative
integers. Let $\LL_{1,n\#}, \ldots, \LL_{n,n\#}, \LL_{\#, n\#}$
denote the line bundles on $\MM_{0,n\#}$ corresponding to the markings
 $1,\ldots,n, \#$ .
Consider the following line bundle on $\MM_{0,n\#}$: 
$$\LL(y_1,\ldots, y_n, y_{\#})\otimes \NN \eqq
\LL^{y_{\#}}_{\#,n\#} \otimes
  \bigotimes_{t=1}^n \LL^{y_t}_{t,n\#} \otimes
 \bigotimes_{i,j} \LL_{i,j}^{x_{i,j}}.$$
where  $\NN=\bigotimes_{i,j} \LL_{i,j}^{x_{i,j}}$ and
$$\LL(y_1,\ldots, y_n, y_{\#})=
\LL^{y_{\#}}_{\#,n\#} \otimes
  \bigotimes_{t=1}^n \LL^{y_t}_{t,n\#}.$$
It suffices to prove for all $k>0$, 
$H^k(\MM_{0,n\#}, \LL(y_1, \ldots, y_n, y_{\#}) \otimes \NN)=0.$

As before, let $\pi_n$ denote the contraction
map $\pi_n: \MM_{0,n\#} \rarr \MM_{0,n}$.
For each $1\leq t \leq n$, 
$$ \LL_{t,n} + D_{t,\#} \eqq \LL_{t,n\#}.$$
The following exact sequences on $\MM_{0,n\#}$ 
therefore exist for
each $1\leq t \leq n$ and non-negative integers $a,b$:
\begin{equation}
\label{redd}
0 \rarr \LL^{a+1}_{t,n}\otimes \LL^b_{t,n\#}
\rarr \LL^{a}_{t,n}\otimes \LL^{b+1}_{t,n\#}
\rarr  \LL^a_{t,n}|_{D_{t,\#}} \rarr 0.
\end{equation}
Via the natural identification $\pi_n:D_{t,\#} \eqq
\M$, the following isomorphisms hold:
$$\LL_{\#, n\#} |_{D_{t,\#}} \eqq \oh,$$
\begin{equation}
\label{rest}
\LL_{t, n\#}|_{D_{t,\#}} \eqq \oh,
\end{equation}
$$\forall s\neq t, \ \LL_{s,n\#} |_{D_{t,\#}} \eqq \LL_{s,n},$$
$$\forall\  1\leq i \leq j, \ 3\leq j \leq n, \ \ 
\LL_{i,j} |_{D_{t,\#}} \eqq \LL_{i,j}.$$
By
repeated use of (\ref{redd}) along with the 
inductive assumption of cohomology
vanishing on $\M \eqq D_{t,\#}$,
it suffices to prove
$$\forall k>0, \ \  H^k(\MM_{0,n\#}, \LL_{\#,n\#}^{y_{\#}}\otimes \NN')=0,$$
$$\NN' \eqq   \bigotimes_{i=1}^n \LL^{y_i}_{i,n} \otimes  \NN.$$
In effect, the sequences (\ref{redd}) are used to
convert  factors of $\LL_{i,n\#}$ to $\LL_{i,n}$
one at a time. The remaining factors of $\LL_{\#,n\#}$
are handled by the following method. 
Since $R^1\pi_{n*}( \LL_{\#,n\#}^{y_{\#}}\otimes \NN')=0$, for all
$k>0$,
$$ H^k(\MM_{0,n\#}, \LL_{\#,n\#}^{y_{\#}}\otimes \NN')\eqq
H^k(\MM_{0,n},\pi_{n*}( \LL_{\#,n\#}^{y_{\#}}\otimes \NN')).$$
Since $\pi_{n*}( \LL_{\#,n\#}^{y_{\#}}\otimes \NN')\eqq
\pi_{n*}( \LL_{\#,n\#}^{y_{\#}})\otimes \NN'$
and $\NN'$ satisfies the vanishing (\ref{f}) by the 
inductive assumption on $\M$, the proof of the
Proposition is completed by Lemma (\ref{ddd}).
\end{pf}
\section{The Proof of Proposition (2)}
The proof of Proposition ($1'$) gives a recursive (in $n$) method
to calculate $$h^0(\M,\ \bigotimes_{i,j} \LL_{i,j}^{x_{i,j}})$$
for non-negative $x_{i,j}$. 
There are 3 reasons why this full recursion is not pursued here:
\begin{enumerate}
\item[(i.)]
The full recursion is complicated.
\item[(ii.)] 
More data is required in the recursion 
than is needed in Proposition (2).
\item[(iii.)] The full recursion
does not respect the symmetry of the variables in
Proposition (2). 
\end{enumerate}
Fortunately there is simple way to sidestep these
problems.

Recall $\g_n(x_1,x_2, \ldots, x_n)=h^0(\M, \LL_{1,n}^{x_1}\otimes
\LL_{2,n}^{x_2} \otimes \cdots \otimes \LL_{n,n}^{x_n})$
is a symmetric function of degree (at most) $n-3$ in
$x_1, \ldots, x_n$. The map
from symmetric functions in $x_1, \ldots, x_n$ to
to symmetric functions in $x_1, \ldots, x_{n-1}$ given by
setting $x_n=0$ is
bijective for symmetric functions of degree at most $n-1$.
Both spaces are spanned by monomials in
$\sigma_1, \ldots, \sigma_{n-1}$ of degree at most $n-1$.
Hence $\g_n$ is completely determined by
$\g_n(x_1, \ldots, x_{n-1}, 0)$.

Suppose the symmetric function $\g_n$ is known.
To determine $\g_{n+1}$, it suffices to
know $\g_{n+1}(x_1, \ldots, x_{n},0)$. The vanishing
of the last entry leads to a great simplification in the
recursion of Proposition ($1'$).
Let 
$$\g_{n+1}(x_1, \ldots, x_n, 0)= h^0(\MM_{0,n\#}, 
\bigotimes_{t=1}^{n}\LL_{t,n\#}^{x_t}).$$
From the sequences (for $1\leq t \leq n$)
\begin{equation}
\label{seqq}
0 \rarr \LL^{a+1}_{t,n}\otimes \LL^b_{t,n\#}
\rarr \LL^{a}_{t,n}\otimes \LL^{b+1}_{t,n\#}
\rarr  \LL^a_{t,n}|_{D_{t,\#}} \rarr 0
\end{equation}
and the restriction equations (\ref{rest}),
the following relation is easily deduced:
$$\g_{n+1}(x_1, \ldots, x_n, 0)=
\g_n(x_1, \ldots, x_n) + 
\sum _{i=1}^{n} \sum_{j=0}^{x_i-1} 
\g_n(x_1, \ldots,x_{i-1}, j, x_{i+1}, 
\ldots, x_{n}).$$
The 
sequences (\ref{seqq}), as before, are used
to convert factors of $\LL_{t,n\#}$ to $\LL_{t,n}$.
By Proposition ($1'$), all higher cohomology
vanishes.
It is the omission of the factor $\LL_{\#,n\#}$ that
simplifies the recursion. 
This concludes the proof of Proposition (2).

\noindent
Department of Mathematics \\
University of Chicago\\
5743 S. University Ave\\
Chicago, IL 60637 \\
rahul@@math.uchicago.edu

\end{document}